\documentclass[11pt]{article}
\usepackage[tbtags]{amsmath}
\usepackage{amssymb}
\usepackage{amsthm}
\usepackage[misc]{ifsym}
\usepackage{cases}
\usepackage{mathrsfs}
\usepackage{graphicx}
\usepackage{subfigure}
\usepackage{float}
\newtheorem{theorem}{Theorem}

\newtheorem{remark}{Remark}

\newtheorem{example}{Example}[section]
\numberwithin{equation}{section}   
\normalsize
\title{\bf A Three-level Stochastic Linear-quadratic Stackelberg Differential Game with Asymmetric Information \thanks{This work is supported by National Natural Science Foundations of China (Grant Nos. 11971266, 11831010, 11571205), and Shandong Provincial Natural Science Foundations (Grant Nos. ZR2020ZD24, ZR2019ZD42).} }

\author{\normalsize  Kaixin Kang\thanks{\it School of Mathematics, Shandong University, Jinan 250100, P.R. China, E-mail: 202012062@mail.sdu.edu.cn} , Jingtao Shi\thanks{\it School of Mathematics, Shandong University, Jinan 250100, P.R. China, E-mail: shijingtao@sdu.edu.cn}}

\begin{document}
\maketitle

\noindent{\bf Abstract:}\quad
This paper is concerned with a three-level stochastic linear-quadratic Stackelberg differential game with asymmetric information, in which three players participate credited as Player 1, Player 2 and Player 3. Player 3 acts as the leader of Player 2 and Player 1, Player 2 acts as the leader of Player 1 and Player 1 acts as the follower. The asymmetric information considered is: the information available to Player 1 is based on the sub-$\sigma$-algebra of the information available to Player 2, and the information available to Player 2 is based on the sub-$\sigma$-algebra of the information available to Player 3. By maximum principle of forward-backward stochastic differential equations and optimal filtering, feedback Stackelberg equilibrium of the game is given with the help of a new system consisting of three Riccati equations.

\vspace{2mm}


\vspace{2mm}


\section{Introduction}

In this paper, we will study a three-level stochastic {\it linear-quadratic} (LQ for short) Stackelberg differential game with asymmetric information. The Stackelberg game, also called leader-follower game, which was proposed by Stackelberg \cite{Stackelberg1952} in 1952, is a kind of game with hierarchical structure. In Stackelberg game, players play the role of the leader or the follower, and make decision sequentially. First of all, we give an example in the supply chain management to introduce the research motivation of this paper.

\begin{example}

(Cooperative advertising and pricing problem) He et al. \cite{HPS2009} studied a cooperative advertising and pricing problem, in which there are two players, a manufacturer and a retailer. Chutani and Sethi \cite{ChutaniSethi2018} considered a cooperative advertising problem under manufacturer and retailer level competition, with a finite number of independent manufacturers and retailers. Kennedy et al. \cite{KSSY2021} extended this problem to the one in a dynamic three-echelon supply chain, which is composed of a manufacturer, a distributor and a retailer. In their supply chain, the manufacturer sells his product to the retailer via the distributor.

We consider the following cooperative advertising and pricing model, which is an extension of that introduced in \cite{KSSY2021}:
\begin{equation}
\left\{
    \begin{aligned}
       dx(t)&=\big[\psi_R(t)\alpha_R(t)+\psi_D(t)\alpha_D(t)+\psi_M(t)\alpha_M(t)\sqrt{1-x(t)}-\delta x(t)\big]dt\\
            &\quad+\sigma(x(t))dW(t),\quad t\in[0,T],\\
        x(0)&=x_0,
    \end{aligned}
\right.
\end{equation}
where $x(\cdot)$ is the market awareness which determines the total sales, the constant $\delta>0$ reflects the rate which potential consumers are lost. $\sigma(x)$ is a variance term, which is usually taken as $\sigma(x)=C_\sigma\sqrt{x(1-x)}$, for some constant $C_\sigma$. The retailer decides the retail price $P_R(\cdot)$ and sets the local advertising effort $\alpha_R(\cdot)$ with the advertising effectiveness $\psi_R(\cdot)$. The distributor decides the distributor price $P_D(\cdot)$ and sets the distributor's advertising effort $\alpha_D(\cdot)$ with the advertising effectiveness $\psi_D(\cdot)$. The manufacturer decides a wholesale price $P_M(\cdot)$, a national advertising effort $\alpha_M(\cdot)$ with the advertising effectiveness $\psi_M(\cdot)$ and a subsidy rate $\phi(\cdot)$ to the retailer's local advertising effort through a vertical cooperative advertising program.

Set $v_R(\cdot)\triangleq(P_R(\cdot), \alpha_R(\cdot), \psi_R(\cdot))$, $v_D(\cdot)\triangleq(P_D(\cdot), \alpha_D(\cdot), \psi_D(\cdot))$ and $v_M(\cdot)\triangleq(P_M(\cdot), \alpha_M(\cdot)$, $\psi_M(\cdot), \phi(\cdot))$, whose values are taken from some admissible control sets $\mathcal{V}_R, \mathcal{V}_D$ and $\mathcal{V}_M$ respectively. Then we encounter a stochastic Stackelberg differential game with three players. In detail, first the retailer's optimal strategy $v^*_R(\cdot)$ is solved by:
\begin{equation}
    \begin{aligned}
J_R(v_R^*(\cdot),v_D(\cdot),v_M(\cdot))=\mathop{max}\limits_{v_R(\cdot)\in\, \mathcal{V}_R} J_R(v_R(\cdot),v_D(\cdot),v_M(\cdot)),\quad \forall\, v_D(\cdot),v_M(\cdot),
    \end{aligned}
\end{equation}
with
\begin{equation}
    \begin{aligned}
&J_R(v_R(\cdot),v_D(\cdot),v_M(\cdot))\\
&=\mathbb{E}\left\{\int_0^T e^{-rt}\big[(P_R(t)-P_D(t))D(P_R(t))x(t)-(1-\phi(t))\alpha_R^2(t)\big]dt\right\}.
    \end{aligned}
\end{equation}
Then the distributor's optimal strategy $v^*_D(\cdot)$ can be obtained by:
\begin{equation}
    \begin{aligned}
J_D(v_R^*(\cdot),v_D^*(\cdot),v_M(\cdot))=\mathop{max}\limits_{v_D(\cdot)\in\, \mathcal{V}_D} J_D(v_R^*(\cdot),v_D(\cdot),v_M(\cdot)),\quad \forall\, v_M(\cdot),
    \end{aligned}
\end{equation}
where
\begin{equation}
    \begin{aligned}
&J_D(v_R(\cdot),v_D(\cdot),v_M(\cdot))\\
&=\mathbb{E}\left\{\int_0^T e^{-rt}\big[(P_D(t)-P_M(t))D(P_R(t))x(t)-k^D D(P_R(t))x(t)-\alpha_D^2(t)\big]dt\right\}.
    \end{aligned}
\end{equation}
Finally, the manufacturer's optimal strategy  $v^*_M$ could be given by
\begin{equation}
    \begin{aligned}
J_M(v_R^*(\cdot),v_D^*(\cdot),v_M^*(\cdot))=\mathop{max}\limits_{v_M(\cdot)\in\, \mathcal{V}_M} J_M(v_R^*(\cdot),v_D^*(\cdot),v_M(\cdot)),
    \end{aligned}
\end{equation}
with
\begin{equation}
    \begin{aligned}
&J_M(v_R(\cdot),v_D(\cdot),v_M(\cdot))\\
&=\mathbb{E}\left\{\int_0^T e^{-rt}\big[(P_M(t)-c)D(P_R(t))x(t)-k^MD(P_R(t))x(t)-\alpha_M^2(t)-\phi(t)\alpha_R^2(t)\big]dt\right\}.
    \end{aligned}
\end{equation}
In the above, $r>0$ is the discount rate, $c>0$ is the manufacturing cost, $k^{M}$ and $k^{D}$ are the transport cost. $0\leq D(p)\leq 1$ is some demand function satisfying usual conditions.

This is a three-level stochastic Stackelberg differential game with three players. Each player hopes to maximize his/her target functional by selecting an appropriate control.
\end{example}

In supply chain management problems, three-level supply chains are often encountered. For example, for a multinational company with multiple sales markets, it is difficult for suppliers to adjust their behavior in direct response to retailers, and the presence of distributors is necessary. This forms a three-level supply chain of suppliers, distributors and retailers.

\begin{figure}[htbp]
  \centering
  \subfigure[]{
  \includegraphics[width=0.15\textwidth]{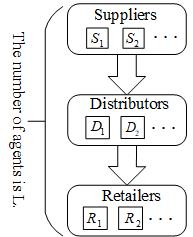}}
  ~~~~~~~~~
  \subfigure[]{
  \includegraphics[width=0.45\textwidth]{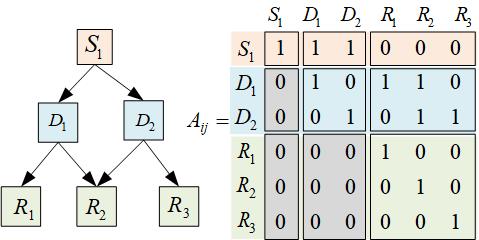}}
  \caption{(a) A schematic of the three-level supply chain; (b)An example of the three-level supply chain}
\end{figure}

A schematic of the three-level supply chain is given in Figure 1.(a), and there can be multiple agents as suppliers, distributors and retails. Another example of the three-level supply chain is given in Figure 1.(b), in which one agent is the supplier, two agents are the distributors, and three agents are the retailers. The matrix can indicate whether there is a leadership relationship between agents, where the leadership relationship means that the information of the ``follower" can be obtained by the ``leader". In Figure 1.(b), for example, $S_{1}$ has a leadership relationship with $D_{1}$, then the position of the matrix $(1,2)$ is 1, $D_{1}$ has no leadership relationship with $R_{3}$, then the position of the matrix $(2,6)$ is 0. Because the leadership of the supply chain is one-way, the shaded part of the matrix must be 0. If no distributors exist, it is a two-level stochastic Stackelberg differential game.

About stochastic Stackelberg differential games with one leader and one follower, there are many related research, such as Yong \cite{Yong2002} and applications to newsvendor-manufacturer problem (\O ksendal et al. \cite{OSU2013}), principal-agent problem (Williams \cite{Williams2015}) and insurer-reinsurer problem (Chen and Shen \cite{CS2018,CS2019}), etc. Mukaidani and Xu \cite{MX2015} studied a stochastic Stackelberg differential game with one leader and multiple followers. Wang and Zhang \cite{wangzhang2020} studied a stochastic LQ Stackelberg differential game of mean-field type with one leader and two followers. Wang and Yan \cite{WangYan2023} researched a Pareto-based stochastic Stackelberg differential game with multi-followers.

However, practically in Stackelberg differential game, due to the emergence of various factors, players often can not observe the complete information, but can only grasp part of the information. This kind of problem is called Stackelberg differential game with {\it asymmetric information}. Shi et al. \cite{SWX2016,SWX2017} studied the two-level stochastic Stackelberg differential games with asymmetric information, in which the information available to the follower is based on the sub-$\sigma$ algebra of that available to the leader. Shi et al. \cite{SWX2020} studied a two-level stochastic LQ Stackelberg differential game with overlapping information, in which the information of the follower and the leader has some overlapping parts, but no mutual inclusion relationship. Li et al. \cite{LMFCM2021} investigated a two-level stochastic LQ Stackelberg differential game under asymmetric information patterns, where the follower uses his observation information to design his strategy whereas the leader implements his strategy using complete information. Zheng and Shi \cite{ZS2022,zhengshi2022} investigated two-level stochastic Stackelberg differential games with partial observation, in which both the leader and the follower have their own observation equations, and the information filtration available to the leader is contained in that to the follower. Yuan et al. \cite{YLH2022} discussed a robust reinsurance contract with asymmetric information in a stochastic Stackelberg differential game. Zhao et al. \cite{ZZGL2022} discussed a stochastic LQ Stackelberg differential game with two leaders and two followers under an incomplete information structure. See more relevant research in the monograph by Ba\c{s}ar and Olsder \cite{BO1999}, and the review paper by Li and Sethi \cite{LiSethi2017}.

Motivated by the above three-level supply chain and the related literatures about the stochastic Stackelberg differential game, in this paper we study a three-level stochastic LQ Stackelberg differential game with asymmetric information. We call the players in the game as Player 1, Player 2 and Player 3. Player 3 acts as the leader of Player 2 and Player 1, Player 2 acts as the leader of Player 1 and Player 1 acts as the follower. The asymmetric information considered is: the information available to Player 1 is based on the sub-$\sigma$-algebra of the information available to Player 2, and the information available to Player 2 is based on the sub-$\sigma$-algebra
of the information available to Player 3. By maximum principle and optimal filtering, feedback Stackelberg equilibrium of the game is given with the help of a new system consisting of three Riccati equations.

The rest of this paper is organized as follows. In Section 2, we formulate our problem. Section 3 is devoted to find the feedback Stackelberg equilibrium of the game. Finally in Section 4, some concluding remarks are given.

\section{Problem formulation}

Let $T>0$ be a finite time horizon, $(\Omega, \mathcal{F}, \mathbb{P})$ be a complete probability space, and $\mathbf{W}=\big\{W_1(t),W_2(t),W_3(t)|t\in[0,T]\big\}$ be a three-dimensional standard Brownian motion defined on it. For $t\in[0,T]$, let $\mathcal{F}_t$ be the natural filtration generated by $\mathbf{W}$ and $\mathcal{F}_t^i\triangleq\sigma\big\{W_i(s);0\leq s\leq t\big\}$ for $i=1,2,3$. In our Stackelberg game, at $t\in[0,T]$ we let $\mathcal{G}_t^1\equiv\mathcal{F}_t^3$ be the information generated by Player 1, $\mathcal{G}_t^2\triangleq\sigma\big\{W_2(s),W_3(s);0\leq s\leq t\big\}$ be the information generated by Player 2, and $\mathcal{G}_t^3\equiv\mathcal{F}_t$ be the information generated by Player 3. Obviously, $\mathcal{G}_t^1\subseteq\mathcal{G}_t^2\subseteq\mathcal{G}_t^3$.

\begin{remark}
Inspired by \cite{ZZGL2022}, we can represent the information owned by the $i$th agent as:
\begin{equation*}
\mathcal{G}_t^i\triangleq\bigcup \limits_{j=1,2,3,\, a_{ij}\neq 0}\mathcal{F}_t^j,\quad \mbox{for }i=1,2,3,
\end{equation*}
where for $i,j=1,2,3$,
\begin{equation*}
a_{ij}\triangleq\left\{\begin{array}{ll}
                       1,&\mathcal{F}_t^j\mbox{ is available to Player i},\\
                       0,&\mathcal{F}_t^j\mbox{ is not available to Player i}.
                       \end{array}
                \right.
\end{equation*}
And the adjacency matrix of the Stackelberg game is
\begin{equation*}
\mathcal{A}=
 \left[
 \begin{array}{ccc}
     0 & 0 & 1 \\
     0 & 1 & 1 \\
     1 & 1 & 1
 \end{array}
 \right].
 \end{equation*}
\end{remark}

We consider the state process $x^{v_1,v_2,v_3}(\cdot)$ which satisfies the following linear {\it stochastic differential equation} (SDE for short):
\begin{equation}\label{asymmstate1}
\left\{
    \begin{aligned}
       dx^{v_1,v_2,v_3}(t)&=\Big[A(t)x^{v_1,v_2,v_3}(t)+\sum_{i=1}^3B_i(t)v_i(t)+b(t)\Big]dt\\
                          &\quad+\sum_{i=1}^3\Big[C_i(t)x^{v_1,v_2,v_3}(t)+\sigma_i(t)\Big]dW_i(t),\quad t\in[0,T],\\
        x^{v_1,v_2,v_3}(0)&=x_0,
    \end{aligned}
\right.
\end{equation}
where $x_0\in\mathbb{R}^n$, $A(\cdot), B_i(\cdot), C_i(\cdot) \in\mathbb{R}^{n\times n}$, $b(\cdot), \sigma_i(\cdot) \in\mathbb{R}^n$ are deterministic and uniformly bounded functions. $v_i(\cdot)\in \mathbb{R}^n$ is the control process of Player $i$, for $i=1,2,3$. The admissible control set $\mathcal{V}_{i}$ of Player $i$ is defined as follows:
\begin{equation}\label{admissible control set}
\mathcal{V}_i\triangleq\left\{v_i:[0,T]\times\Omega\rightarrow \mathbb{R}^n|\, v_i(\cdot)\mbox{ is }\mathcal{G}_t^i\mbox{-adapted, and }\mathbb{E}\int_0^T|v_i(t)|^2dt<\infty\right\},\,\, i=1,2,3.
\end{equation}

The quadratic cost functional of $i$th Player $(i=1,2,3)$ is defined as follows:
\begin{equation}\label{cost fuinctional}
\begin{split}
&J_i(v_1(\cdot),v_2(\cdot),v_3(\cdot))\triangleq\frac{1}{2}\mathbb{E}\bigg\{\int_0^T\Big[\big\langle Q_i(t)x^{v_1,v_2,v_3}(t),x^{v_1,v_2,v_3}(t)\big\rangle+\big\langle R_i(t)v_i(t),v_i(t)\big\rangle\\
&\qquad +2\big\langle m_i(t),x^{v_1,v_2,v_3}(t)\big\rangle+2\big\langle n_i(t),v_i(t)\big\rangle\Big]dt+\big\langle G_i(T)x^{v_1,v_2,v_3}(T),x^{v_1,v_2,v_3}(T)\big\rangle \bigg\},
\end{split}
\end{equation}
where $Q_i(\cdot), R_i(\cdot), G_i\in \mathbb{R}^{n\times n}, m_i(\cdot), n_i(\cdot)\in \mathbb{R}^n\,\, (i=1,2,3)$ are deterministic and uniformly bounded functions. In addition, $R_i(\cdot)\gg0, Q_i(\cdot)\geq0, G_i\geq0$, for $i=1,2,3.$

The definition of a three-level Stackelberg game's equilibrium strategy is given in \cite{BO1999}, Definition 3.37. In detail, the players' optimal goals are as follows.

First, for Player 1:
\begin{equation}
J_1(v_1^*(\cdot),v_2(\cdot),v_3(\cdot))=\inf\limits_{v_1(\cdot)\in\, \mathcal{V}_1}J_1(v_1(\cdot),v_2(\cdot),v_3(\cdot)),\ \forall v_2(\cdot)\in\mathcal{V}_2,\,\, v_3(\cdot)\in\mathcal{V}_3;
\end{equation}
then, for Player 2:
\begin{equation}
J_2(v_1^*(\cdot),v_2^*(\cdot),v_3(\cdot))=\inf\limits_{v_2(\cdot)\in\, \mathcal{V}_2}J_2(v_1^*(\cdot),v_2(\cdot),v_3(\cdot)),\ \forall v_3(\cdot)\in\mathcal{V}_3;
\end{equation}
finally, for Player 3:
\begin{equation}
J_3(v_1^*(\cdot),v_2^*(\cdot),v_3^*(\cdot))=\inf\limits_{v_3(\cdot)\in\, \mathcal{V}_3}J_3(v_1^*(\cdot),v_2^*(\cdot),v_3(\cdot)).
\end{equation}

If such $(v_1^*(\cdot), v_2^*(\cdot), v_3^*(\cdot))$ exists, it is called an {\it open-loop Stackelberg equilibrium strategy} of our stochastic LQ Stackelberg differential game with asymmetric infomation. In this paper, our ultimate goal is to find the optimal feedback strategies of Player 1, Player 2 and Player 3.

\section{Stackelberg equilibrium strategies}

In this section, the problem is resolved in three steps, to solve optimal strategies of Player 1, Player 2 and Player 3 in turn. The main tool is various stochastic versions of {\it Pontryagin's maximum principle}. For any $\mathcal{F}_t$-adapted process $\xi(\cdot)$, we denote by
\begin{equation}\label{filter estimates}
\hat{\xi}(t)\triangleq\mathbb{E}[\xi(t)|\mathcal{G}_t^2],\quad \check{\xi}(t)\triangleq\mathbb{E}[\xi(t)|\mathcal{G}_t^1]
\end{equation}
its optimal filtering estimates, for $t\in[0,T]$.

\subsection{Problem of Player 1}

First, we define the Hamiltionian function of Player 1 $H_1:[0,T]\times\mathbb{R}^n\times\mathbb{R}^n\times\mathbb{R}^n\times\mathbb{R}^n\times\mathbb{R}^n\times\mathbb{R}^n\times\mathbb{R}^n\times\mathbb{R}^n\rightarrow\mathbb{R}$ as
\begin{equation}\label{Hamiltonian player 1}
\begin{split}
  &H_1(t,x,v_1,v_2,v_3,y_1,z_1,z_2,z_3)\\
  &\triangleq\left[A(t)x+\sum_{i=1}^3 B_i(t)v_i+b(t)\right]^\top y_1+\sum_{i=1}^3 \big[C_i(t)x+\sigma_i(t)\big]^\top z_i\\
  &\quad-\bigg[\frac{1}{2}\langle Q_1(t)x,x\rangle+\frac{1}{2}\langle R_1(t)v_1,v_1\rangle+\langle m_1(t),x\rangle+\langle n_1(t),v_1\rangle\bigg].
\end{split}
\end{equation}
Thus, the optimal strategy of Player 1 can be obtained by the stochastic maximum principle with partial information (See, for example, Baghery and \O ksendal \cite{AO2007}, or Shi et al. \cite{SWX2016}):
\begin{equation}\label{asymmv1chu}
 R_1(t)v_1^*(t)-B_1(t)^\top\check{y}_1(t)+n_1(t)=0,\quad a.e.\ t\in[0,T],\ \ a.s.,
\end{equation}
where the $\mathcal{F}_{t}$-adapted process quadruple ($y_1(\cdot), z_1(\cdot), z_2(\cdot), z_3(\cdot)$) is the unique solution to the {\it backward SDE} (BSDE for short):
\begin{equation}\label{asymmad1}
\left\{
    \begin{aligned}
       -dy_1(t)&=\left[A(t)^\top y_1(t)+\sum_{i=1}^3 C_i(t)^\top z_i(t)-Q_1(t)x_1^*(t)-m_1(t)\right]dt\\
               &\quad -\sum_{i=1}^3 z_i(t)dW_i(t),\quad t\in[0,T],\\
         y_1(T)&=-G_1x_1^*(T),
    \end{aligned}
\right.
\end{equation}
where for simplicity, we denote $x^*_1(\cdot)\equiv x^{v_1^*,v_2,v_3}(\cdot)$.

We want to represent the optimal strategy in the form of state feedback. For this target, let
\begin{equation}\label{asymmy1}
y_1(t)=-p(t)x_1^*(t)-\varphi(t),\quad t\in[0,T],
\end{equation}
where $p(\cdot)$ is a deterministic, differentiable function with $p(T)=G_1$, and $\mathcal{F}_t$-adapted processes pair $(\varphi(\cdot), \theta(\cdot) )$ satisfies the BSDE:
\begin{equation}\label{variphi}
\left\{
    \begin{aligned}
       d\varphi(t)&=\alpha(t)dt+\theta(t)dW_3(t),\quad t\in[0,T],\\
        \varphi(T)&=0,
    \end{aligned}
\right.
\end{equation}
where $\alpha(\cdot)$ is some $\mathcal{F}_t$-adapted process to be determined.

Let (\ref{asymmy1}) be filtered on $\mathcal{G}_{t}^{1}$ to obtain $\check{y}_1(\cdot)$. Substituting it into (\ref{asymmv1chu}), we can get the optimal state feedback strategy of Player 1 as:
\begin{equation}\label{asymmvbiaoda1}
v_1^*(t)=-R_1^{-1}(t)\left[B_1(t)^\top p(t)\check{x}_1^*(t)+B_1(t)^\top\check{\varphi}(t)+n_{1}(t)\right],\quad a.e.\ t\in[0,T],\ \ a.s..
\end{equation}

In the same way as Shi et al. \cite{SWX2017}, we apply It\^o's formula to (\ref{asymmy1}). Then, making a comparison between that and (\ref{asymmad1}), we obtain the following Riccati equation which admits a unique differentiable solution $p(\cdot)$ (See, for example, Yong and Zhou \cite{YZ1999}, Chapter 6):
\begin{equation}\label{asymmRiccati1}
\left\{
    \begin{aligned}
      &\dot{p}(t)+p(t)A(t)+A(t)^\top p(t)+\sum_{i=1}^3 C_i(t)^\top p(t)C_i(t)\\
      &\quad -p(t)B_1(t)R_1^{-1}(t)B_1(t)^\top p(t)+Q_1(t)=0,\quad t\in[0,T],\\
      &p(T)=G_1.
    \end{aligned}
\right.
\end{equation}
Moreover, we can get that $(\varphi(\cdot),\theta(\cdot))$ satisfies the following BSDE:
\begin{equation}\label{asymmvarphi}
\left\{
    \begin{aligned}
     -d\varphi(t)&=\bigg[p(t)B_1(t)R_1^{-1}(t)B_1(t)^\top p(t)x_1^*(t)-p(t)B_1(t)R_1^{-1}(t)B_1(t)^\top p(t)\check{x}_1^*(t)\\
                 &\qquad +A(t)^\top\varphi(t)-p(t)B_1(t)R_1^{-1}(t)B_1(t)^\top\check{\varphi}(t)+C(t)^\top\theta(t)\\
                 &\qquad +\sum_{i=2}^3 p(t)B_i(t)v_i(t)+\sum_{i=1}^3 C_i(t)^\top p(t)\sigma_i(t)+p(t)b(t)+m_1(t)\\
                 &\qquad -p(t)B_1(t)R_1^{-1}(t)n_1(t)\bigg]dt-\theta(t)dW_3(t),\quad t\in[0,T],\\
       \varphi(T)&=0.
    \end{aligned}
\right.
\end{equation}

\begin{remark}
Due to the length of the article, details of the calculation in (\ref{asymmRiccati1}), (\ref{asymmvarphi}) are omitted. Readers can refer to \cite{SWX2017}. At this point, it is equivalent to decouple $x^*(\cdot)$ and $y_1(\cdot)$, which ensures that (\ref{asymmstate1}), together with (\ref{asymmad1}), admits a unique solution $(x_1^*(\cdot),y_1(\cdot),z_1(\cdot),z_2(\cdot),z_3(\cdot))$.
\end{remark}

Put the optimal strategy $v_1^*(\cdot)$ of Player 1 into the state equation (\ref{asymmstate1}) of $x_1^*(\cdot)$, and set
\begin{equation}
    \begin{aligned}
      \bar{A}(t)&\triangleq A(t)-B_1(t)R_1^{-1}(t)B_1(t)^\top p(t),\quad \bar{F}_{1}(t)\triangleq-B_1(t)R_1^{-1}(t)B_1(t)^\top,\\
      \bar{F}_{2}(t)&\triangleq B_2(t)^\top p(t),\quad \bar{F}_3(t)\triangleq B_3(t)^\top p(t),\quad \bar{b}(t)\triangleq b(t)-B_1(t)R_1^{-1}(t)n_1(t),\\
      \bar{f}_{1}(t)&\triangleq p(t)b(t)+\sum_{i=1}^3 C_i(t)^\top p(t)\sigma_i(t)+m_1(t)-p(t)B_1(t)R_1^{-1}(t)n_1(t),\quad t\in[0,T].
    \end{aligned}
\end{equation}
Applying Lemma 5.4 in Xiong \cite{Xiong2008} to state equation of $x_1^*(\cdot)$ and (\ref{asymmvarphi}), we derive the optimal filtering equation:
\begin{equation}\label{asymmwenti1}
\left\{
    \begin{aligned}
      d\check{x}_1^*(t)&=\left[\bar{A}(t)\check{x}_1^*(t)+\bar{F}_1(t)\check{\varphi}(t)+\sum_{i=2}^3 B_i(t)\check{v}_i(t)+\bar{b}(t)\right]dt+\big[C_3(t)\check{x}_1^*(t)+\sigma_3(t)\big]dW_3(t),\\
   -d\check{\varphi}(t)&=\left[\bar{A}(t)^\top\check{\varphi}(t)+C_3(t)^\top\check{\theta}(t)+\sum_{i=2}^3\bar{F}_i(t)^\top\check{v}_i(t)+\bar{f}_1(t)\right]dt-\check{\theta}(t)dW_3(t),\quad t\in[0,T],\\
       \check{x}_1^*(0)&=x_0,\quad \check{\varphi}(T)=0,
    \end{aligned}
\right.
\end{equation}
which, obviously, admits a unique $\mathcal{G}_t^1$-adapted solution $(\check{x}_1^*(\cdot), \check{\varphi}(\cdot), \check{\theta}(\cdot))$.

For any $v_2(\cdot)\in\mathcal{V}_2$ and $v_3(\cdot)\in\mathcal{V}_3$, the problem of Player 1, first, is solved in the following theorem.
\begin{theorem}
Let $p(\cdot)$ satisfy (\ref{asymmRiccati1}), Player 1's optimal strategy $v_1^*(\cdot)$ is given by (\ref{asymmvbiaoda1}), where $(\check{x}_1^*(\cdot), \check{\varphi}(\cdot), \check{\theta}(\cdot))$ is the unique $\mathcal{G}_t^1$-adapted solution to (\ref{asymmwenti1}).
\end{theorem}

\subsection{Problem of Player 2}

This section focuses on solving the optimal control problem of Player 2. After Player 1 exercises the optimal strategy $v_1^*(\cdot)=v_1^*[\cdot;v_2(\cdot),v_3(\cdot)]$ (which means that $v_1^*(\cdot)$ depends on $v_2(\cdot)$ and $v_3(\cdot)$) that we got in the previous subsection, the state equation of Player 2's problem becomes:
\begin{equation}\label{asymmstate2}
\left\{
    \begin{aligned}
             dx_2(t)&=\left[A(t)x_2(t)+\big(\bar{A}(t)-A(t)\big)\check{x}_2(t)+\bar{F}_1(t)\check{\varphi}(t)+\sum_{i=2}^3 B_i(t)v_i(t)+\bar{b}(t)\right]dt\\
                    &\quad +\sum_{i=1}^3 \big[C_i(t)x_2(t)+\sigma_i(t)\big]dW_i(t),\\
-d\check{\varphi}(t)&=\left[\bar{A}(t)^\top\check{\varphi}(t)+C_3(t)^\top\check{\theta}(t)+\sum_{i=2}^3\bar{F}_i(t)^\top\check{v}_i(t)+\bar{f}_1(t)\right]dt-\check{\theta}(t)dW_3(t),\quad t\in[0,T],\\
              x_2(0)&=x_0,\quad \check{\varphi}(T)=0,
    \end{aligned}
\right.
\end{equation}
where $x_2(\cdot)\triangleq x^{v_1^*,v_2,v_3}(\cdot)(\equiv x_1^*(\cdot)), \check{x}_2(\cdot)\triangleq\check{x}^{v_1^*,v_2,v_3}(\cdot)$. Noting that equation (\ref{asymmstate2}) is decoupled {\it conditional mean-field forward-backward SDE} (CMFFBSDE for short) (see, for example, \cite{SWX2016}), which admits a unique adapted solution $(x_2(\cdot),\check{\varphi}(\cdot),\check{\theta}(\cdot))$.

Define the Hamiltionian function $H_2:[0,T]\times\mathbb{R}^n\times\mathbb{R}^n\times\mathbb{R}^n\times\mathbb{R}^n\times\mathbb{R}^n\times\mathbb{R}^n\times\mathbb{R}^n\times\mathbb{R}^n\times\mathbb{R}^n\times\mathbb{R}^n\rightarrow\mathbb{R}$ of Player 2 as
\begin{equation}\label{Hamiltonian player 2}
\begin{aligned}
  &H_2(t,x_2,\check{\varphi},\check{\theta},v_2,v_3,y_2,\tilde{z}_1,\tilde{z}_2,\tilde{z}_3,\psi_2)\\
  &\triangleq\left[A(t)x_2+(\bar{A}(t)-A(t))\check{x}_2+\bar{F}_1(t)\check{\varphi}+\sum_{i=2}^3 B_i(t)v_i+\bar{b}(t)\right]^{\top}y_2\\
  &\quad +\sum_{i=1}^3 \big[C_i(t)x_2+\sigma_i(t)\big]^\top \tilde{z}_i(t)+\left[\bar{A}(t)^\top\check{\varphi}+C_3(t)^\top\check{\theta}+\sum_{i=2}^3 \bar{F}_i(t)^\top\check{v}_i+\bar{f}_1(t)\right]^\top\psi_2\\
  &\quad +\left[\frac{1}{2}\langle Q_2(t)x_2,x_2\rangle+\frac{1}{2}\langle R_2(t)v_2,v_2\rangle+\langle m_2(t),x_2\rangle+\langle n_2(t),v_2\rangle\right].
\end{aligned}
\end{equation}
The optimal strategy $v_2^*(\cdot)$ of Player 2 can be obtained by the maximum principle of FBSDE with asymmetric information (see for example, \cite{SWX2016}, \cite{SWX2017}):
\begin{equation}\label{}
  R_2(t)v_2^*(t)+B_2(t)^\top\hat{y}_2(t)+\bar{F}_2(t)\check{\psi}_2(t)+n_2(t)=0,\quad a.e.\ t\in[0,T],\ \ a.s.,
\end{equation}
where $(y_2(\cdot), \tilde{z}_1(\cdot), \tilde{z}_2(\cdot), \tilde{z}_3(\cdot), \psi_2(\cdot))$ is the unique $\mathcal{F}_t$-adapted solution to the forward FBSDE:
\begin{equation}\label{asymmad2}
\left\{
    \begin{aligned}
      -dy_2(t)&=\left[A(t)^\top y_2(t)+(\bar{A}(t)-A(t))^\top\check{y}_2(t)+\sum_{i=1}^3 C_i(t)^\top\tilde{z}_i(t)\right.\\
              &\qquad +Q_2(t)x_2^*(t)+m_2(t)\bigg]dt-\sum_{i=1}^3\tilde{z}_i(t)dW_i(t),\\
    d\psi_2(t)&=\big[\bar{A}(t)\psi_2(t)+\bar{F}_1(t)y_2(t)\big]dt+C_3(t)\psi_2(t)dW_3(t),\quad t\in[0,T],\\
        y_2(T)&=G_2x_2^*(T),\quad \psi_2(0)=0,
    \end{aligned}
\right.
\end{equation}
with $(x_2^*(\cdot)\triangleq x^{v_1^*,v_2^*,v_3}(\cdot), \check{\varphi}^*(\cdot), \check{\theta}^*(\cdot))$ is the solution to (\ref{asymmstate2}) corresponding to $v_2^*(\cdot)$.

Noting that (\ref{asymmad2}) is a coupled CMFFBSDE, we need to decouple it to ensure its solvability.
Inspired by Yong \cite{Yong2002}, let $(x_2^*(\cdot),\psi_2(\cdot))^\top$ be the optimal ``state", and set
\begin{equation}
    \begin{aligned}
    &X_2\triangleq\left(\begin{matrix}x_2^*\\\psi_2\end{matrix}\right),\
     Y_2\triangleq\left(\begin{matrix}y_2\\\check{\varphi}^*\end{matrix}\right),\
     \tilde{Z}_1\triangleq\left(\begin{matrix}\tilde{z}_1\\0\end{matrix}\right),\
     \tilde{Z}_2\triangleq\left(\begin{matrix}\tilde{z}_2\\0\end{matrix}\right),\
     \tilde{Z}_3\triangleq\left(\begin{matrix}\tilde{z}_3\\\check{\theta}^*\end{matrix}\right),\\
    &X_0\triangleq\left(\begin{matrix}x_0\\0\end{matrix}\right),\
     \bar{\mathcal{A}}_1\triangleq\left(\begin{matrix}A& 0\\0 &\bar{A}\end{matrix}\right),\
     \bar{\mathcal{A}}_2\triangleq\left(\begin{matrix}\bar{A}-A& 0\\0 & 0\end{matrix}\right),\
     \bar{\mathcal{F}}_1\triangleq\left(\begin{matrix}0& \bar{F}_1\\\bar{F}_1& 0\end{matrix}\right),\\
    &\bar{\mathcal{B}}_2\triangleq\left(\begin{matrix}B_2\\0\end{matrix}\right),\
     \bar{\mathcal{B}}_3\triangleq\left(\begin{matrix}B_3\\0\end{matrix}\right),\
     \bar{\mathcal{C}}_1\triangleq\left(\begin{matrix}C_1& 0\\0& 0\end{matrix}\right),\
     \bar{\mathcal{C}}_2\triangleq\left(\begin{matrix}C_2& 0\\0& 0\end{matrix}\right),\
     \bar{\mathcal{C}}_3\triangleq\left(\begin{matrix}C_3& 0\\0& C_3\end{matrix}\right),\\
    &\bar{\mathcal{Q}}_2\triangleq\left(\begin{matrix}Q_2& 0\\0& 0\end{matrix}\right),\
     \bar{\mathcal{G}}_2\triangleq\left(\begin{matrix}G_2& 0\\ 0& 0\end{matrix}\right),\
     \bar{\mathcal{F}}_2\triangleq\left(\begin{matrix}0&\bar{F}_2\end{matrix}\right),\
     \bar{\mathcal{F}}_3\triangleq\left(\begin{matrix} 0&\bar{F}_3\end{matrix}\right),\
     \bar{b}_2=\left(\begin{matrix}\bar{b}\\0\end{matrix}\right),\\
    &\bar{\sigma}_1\triangleq\left(\begin{matrix}\sigma_1\\0\end{matrix}\right),\
     \bar{\sigma}_2\triangleq\left(\begin{matrix}\sigma_2\\0\end{matrix}\right),\
     \bar{\sigma}_3\triangleq\left(\begin{matrix}\sigma_3\\0\end{matrix}\right),\
     \bar{f}_2\triangleq\left(\begin{matrix}m_2\\\bar{f}_1\end{matrix}\right).
    \end{aligned}
\end{equation}
Noting that time variables are usually omitted somewhere since now, for notational simplicity. Then (\ref{asymmstate2}) (corresponding to $v_2^*(\cdot)$) and (\ref{asymmad2}) can be rewritten as:
\begin{equation}\label{asymmkuowei2}
\left\{
    \begin{aligned}
       dX_2(t)&=\big[\bar{\mathcal{A}}_1X_2+\bar{\mathcal{A}}_2\check{X}_2+\bar{\mathcal{F}}_1Y_2+\bar{\mathcal{B}}_2v_2^*+\bar{\mathcal{B}}_3v_3+\bar{b}_2\big]dt
               +\sum_{i=1}^3\big[\bar{\mathcal{C}}_iX_2+\bar{\sigma}_i\big]dW_i(t),\\
      -dY_2(t)&=\Big[\bar{\mathcal{Q}}_2X_2+\hat{\mathcal{A}}_1^\top Y_2+\hat{\mathcal{A}}_2^\top\check{Y}_2+\sum_{i=1}^3 \bar{\mathcal{C}}_i^\top Z_i
               +\bar{\mathcal{F}}_2^\top\check{v}_2^*+\bar{\mathcal{F}}_3^\top\check{v}_3+\bar{f}_2\Big]dt\\
              &\quad -\sum_{i=1}^3 \tilde{Z}_i(t)dW_i(t),\quad t\in[0,T],\\
        X_2(0)&=X_0,\quad Y_2(T)=\mathcal{G}_2X_2(T),
    \end{aligned}
\right.
\end{equation}
and similarly, the optimal strategy $v_2^*(\cdot)$ of Player 2 can be written as:
\begin{equation}\label{asymmv2chu}
  v_2^*(t)=-R_2^{-1}\big[\bar{\mathcal{B}}_2^\top\hat{Y}_2+\bar{\mathcal{F}}_2\check{X}_2+n_2\big],\quad a.e.\ t\in[0,T],\ \ a.s..
\end{equation}

Next, in order to represent the optimal strategy $v_2^*(\cdot)$ of Player 2 in the form of state feedback, we assume
\begin{equation}\label{asymmy2}
Y_2(t)=P_1(t)X_2(t)+P_2(t)\check{X}_2(t)+\Phi(t),\quad t\in[0,T],
\end{equation}
where $P_1(\cdot), P_2(\cdot)$ are deterministic, differentiable functions with $P_1(T)=G_2, P_2(T)=0$, and $\mathcal{F}_t$-adapted processes pair $(\Phi(\cdot),\Lambda(\cdot))$ satisfies the BSDE:
\begin{equation}
\left\{
    \begin{aligned}
       d\Phi(t)&=\Gamma(t)dt+\Lambda(t)dW_2(t),\\
        \Phi(T)&=0,
    \end{aligned}
\right.
\end{equation}
where $\Gamma(\cdot)$ is some $\mathcal{F}_t$-adapted process to be determined later.

Let (\ref{asymmy2}) be filtered on $\mathcal{G}_t^2$ to obtain $\hat{Y}_2(\cdot)$. Substituting it into (\ref{asymmv2chu}), we can represent the optimal state feedback strategy $v_2^*(\cdot)$ of Player 2 as:
\begin{equation}\label{asymmvbiaoda21}
\begin{aligned}
v_2^*(t)&=-R_2^{-1}\left\{\bar{\mathcal{B}}_2^\top P_1\hat{X}_2+(\bar{\mathcal{B}}_2^\top P_2+\bar{\mathcal{F}}_2)\check{X}_2+\bar{\mathcal{B}}_2^\top\hat{\Phi}+n_2\right\},\quad a.e.\ t\in[0,T],\ \ a.s..
\end{aligned}
\end{equation}

Using the same method as in Section 3.1, we introduce a system of two Riccati equations as follows:
\begin{equation}\label{asymmRiccati21}
\left\{
    \begin{aligned}
      &\dot{P}_1+P_1\bar{\mathcal{A}}_1+\bar{\mathcal{A}}_1^\top P_1+\sum_{i=1}^3 \bar{\mathcal{C}}_i^\top P_1\bar{\mathcal{C}}_i-\bar{\mathcal{F}}_2^\top R_2^{-1}\bar{\mathcal{B}}_2^\top P_1
       +P_1\big(\bar{\mathcal{F}}_1-\bar{\mathcal{B}}_2R_2^{-1}\bar{\mathcal{B}}_2^\top\big)P_1+\bar{\mathcal{Q}}_2=0,\\
      &P_1(T)=\mathcal{G}_2,
    \end{aligned}
\right.
\end{equation}
\begin{equation}\label{asymmRiccati22}
\left\{
    \begin{aligned}
      &\dot{P}_2+\big(\bar{\mathcal{A}}_1^\top+\bar{\mathcal{A}}_2^\top-\bar{\mathcal{F}}_2^\top R_2^{-1}\bar{\mathcal{B}}_2^\top\big)P_2
       +P_2\big(\bar{\mathcal{A}}_1+\bar{\mathcal{A}}_2-\bar{\mathcal{B}}_2R_2^{-1}\bar{\mathcal{F}}_2^\top\big)\\
      &\quad +P_1\big(\bar{\mathcal{F}}_1-\bar{\mathcal{B}}_2R_2^{-1}\bar{\mathcal{B}}_2^\top\big)P_2
       +P_2\big(\bar{\mathcal{F}}_1-\bar{\mathcal{B}}_2R_2^{-1}\bar{\mathcal{B}}_2^\top\big)P_2+\bar{\mathcal{C}}_2^\top P_2\bar{\mathcal{C}}_3\\
      &\quad +P_1\big(\bar{\mathcal{A}}_2-\bar{\mathcal{B}}_2R_2^{-1})\bar{\mathcal{F}}_2^\top\big)-\bar{\mathcal{F}}_2^\top R_2^{-1}\bar{F}_2^\top=0,\\
      &P_2(T)=0.
    \end{aligned}
\right.
\end{equation}
Noting system (\ref{asymmRiccati21}), (\ref{asymmRiccati22}) is not coupled, it admits a differentiable solution pair $(P_1(\cdot), P_2(\cdot))$. Then, we can get $(\Phi(\cdot),\Lambda(\cdot))$ satisfies the following BSDE:
\begin{equation}\label{asymmPhi21}
\left\{
    \begin{aligned}
      -d\Phi(t)&=\Big\{\big[P_1\bar{\mathcal{B}}_2R_2^{-1}\bar{\mathcal{B}}_2^\top P_1+\bar{\mathcal{F}}_2^\top R_2^{-1}\bar{\mathcal{B}}_2^\top P_1\big]X_2
                -\big[P_1\bar{\mathcal{B}}_2R_2^{-1}\bar{\mathcal{B}}_2^\top P_1+\bar{\mathcal{F}}_2^\top R_2^{-1}\bar{\mathcal{B}}_2^\top P_1\big]\hat{X}_2\\
               &\qquad +\big[\bar{\mathcal{A}}_1^\top+P_1\bar{\mathcal{F}}_1\big]\Phi
                -\big[P_1\bar{\mathcal{B}}_2R_2^{-1}\bar{\mathcal{B}}_2^\top+\bar{\mathcal{F}}_2^\top R_2^{-1}\bar{\mathcal{B}}_2^\top\big]\hat{\Phi}\\
               &\qquad +\big[\bar{\mathcal{A}}_2^\top+P_2\big(\bar{\mathcal{F}}_1-\bar{\mathcal{B}}_2R_2^{-1}\bar{\mathcal{B}}_2^\top\big)\big]\check{\Phi}
                +\big[P_1\bar{\mathcal{B}}_3+\bar{\mathcal{F}}^\top_3\big]v_3+P_2\bar{\mathcal{B}}_3\check{v}_3+\bar{\mathcal{C}}_2^\top\Lambda\\
               &\qquad +\big(P_1+P_2)\bar{b}_2+\bar{\mathcal{C}}_1^\top P_1\bar{\sigma}_1+\bar{\mathcal{C}}_2^\top P_1\bar{\sigma}_2+\bar{\mathcal{C}}_3^\top\big(P_1+P_2\big)\bar{\sigma}_3\\
               &\qquad -\big(P_1+P_2\big)\bar{\mathcal{B}}_2R_2^{-1}n_2+\bar{f}_2-\bar{\mathcal{F}}_2^\top R_2^{-1}n_2\Big\}dt-\Lambda(t)dW_2(t),\quad t\in[0,T],\\
       \Phi(T)&=0.
    \end{aligned}
\right.
\end{equation}

\begin{remark}
Readers can refer to \cite{SWX2017} for more details about derivation of (\ref{asymmRiccati21}) and (\ref{asymmRiccati22}). In this way it can also be shown the solvability of (\ref{asymmkuowei2}), since $(\hat{X}_2(\cdot), \check{X}_2(\cdot), \hat{\Phi}(\cdot),\hat{\Lambda}(\cdot))$ in (\ref{asymmvbiaoda21}) are unique $\mathcal{G}_t^2$-adapted solutions to the following filter equations (\ref{asymmwenti21}), (\ref{asymmwenti22}). Thus, the solvability of CMFFBSDE (\ref{asymmad2}) can be obtained.
\end{remark}

In fact, putting the optimal strategy $v_2^*(\cdot)$ of Player 2 (\ref{asymmvbiaoda21}) into the equation (\ref{asymmkuowei2}) of the state $X_2(\cdot)$, setting
\begin{equation*}
    \begin{aligned}
      \bar{\bar{\mathcal{A}}}_1&\triangleq\bar{\mathcal{A}}_1+\bar{\mathcal{F}}_1P_1,\quad \bar{\bar{\mathcal{A}}}_2\triangleq-\bar{\mathcal{B}}_2R_2^{-1}\bar{\mathcal{B}}_2^\top P_1,\
      \bar{\bar{\mathcal{A}}}_3\triangleq\bar{\mathcal{A}}_2+\bar{\mathcal{F}}_1P_2-\bar{\mathcal{B}}_2R_2^{-1}\bar{\mathcal{B}}_2^\top P_2-\bar{\mathcal{B}}_2R_2^{-1}\bar{\mathcal{F}}_2,\\
      \bar{\bar{\mathcal{F}}}_1&\triangleq\bar{\mathcal{F}}_1-\bar{\mathcal{B}}_2R_2^{-1}\bar{\mathcal{B}}_2^\top,\quad \bar{\bar{\mathcal{F}}}_3\triangleq P_1\bar{\mathcal{B}}_3+\bar{\mathcal{F}}_3^\top,\quad
                                \bar{\bar{b}}_2\triangleq\bar{b}_2-\bar{\mathcal{B}}_2R_2^{-1}n_2,\\
                \bar{\bar{f}}_2&\triangleq\bar{f}_2+(P_1+P_2)\bar{b}_2+\bar{\mathcal{C}}_1^\top P_1\bar{\sigma}_1
                                +\bar{\mathcal{C}}_2^\top P_1\bar{\sigma}_2+\bar{\mathcal{C}}_3^\top(P_1+P_2)\bar{\sigma}_3\\
                               &\quad -(P_1+P_2)\bar{\mathcal{B}}_2R_2^{-1}n_2-\bar{\mathcal{F}}_2^\top R_2^{-1}n_2,\\
                    \mathcal{H}&\triangleq P_1\bar{\mathcal{B}}_2R_2^{-1}\bar{\mathcal{B}}_2\top P_1+\bar{\mathcal{F}}_2^\top R_2^{-1}\bar{\mathcal{B}}_2^\top P_1,
    \end{aligned}
\end{equation*}
and applying Lemma 5.4 in \cite{Xiong2008} to it and (\ref{asymmPhi21}), we derive the optimal filtering equations:
\begin{equation}\label{asymmwenti21}
\left\{
    \begin{aligned}
  d\hat{X}_2(t)&=\Big\{\big(\bar{\bar{\mathcal{A}}}_1+\bar{\bar{\mathcal{A}}}_2\big)\hat{X}_2+\bar{\bar{\mathcal{A}}}_3\check{X}_2+\bar{\bar{\mathcal{F}}}_1\hat{\Phi}
                +\bar{\mathcal{B}}_3\hat{v}_3+\bar{\bar{b}}_2\Big\}dt +\sum_{i=2}^3 \big[\bar{\mathcal{C}}_i\hat{X}_2+\bar{\sigma}_i\big]dW_i(t),\\
-d\hat{\Phi}(t)&=\Big\{\big[\bar{\bar{\mathcal{A}}}_1^\top+\bar{\bar{\mathcal{A}}}_2^\top-\bar{\mathcal{F}}_2^\top R_2^{-1}\bar{\mathcal{B}}_2^\top\big]\hat{\Phi}
                +\big[\bar{\bar{\mathcal{A}}}_3^\top+\bar{\mathcal{F}}_2^\top R_2^{-1}\bar{\mathcal{B}}_2^\top]\check{\Phi}
                +\bar{\mathcal{C}}_2^\top\hat{\Lambda}+\bar{\bar{\mathcal{F}}}_3^\top\hat{v}_3\\
               &\qquad +\big(\bar{\bar{\mathcal{F}}}_3^\top-\bar{\mathcal{F}}_3^\top\big)\check{v}_3+\bar{\bar{f}}_2\Big\}dt-\hat{\Lambda}(t)dW_2(t),\quad t\in[0,T],\\
 \hat{X}_{2}(0)&=X_0,\quad \hat{\Phi}_2(T)=0,
    \end{aligned}
\right.
\end{equation}
and
\begin{equation}\label{asymmwenti22}
\left\{
    \begin{aligned}
  d\check{X}_2(t)&=\Big\{\big(\bar{\bar{\mathcal{A}}}_1+\bar{\bar{\mathcal{A}}}_2+\bar{\bar{\mathcal{A}}}_3\big)\check{X}_2+\bar{\bar{\mathcal{F}}}_1\check{\Phi}
                  +\bar{\mathcal{B}}_3\check{v}_3+\bar{\bar{b}}_2\Big\}dt+\big[\bar{\mathcal{C}}_3\check{X}_2+\bar{\sigma}_3\big]dW_3(t),\\
-d\check{\Phi}(t)&=\Big\{\big(\bar{\bar{\mathcal{A}}}_1+\bar{\bar{\mathcal{A}}}_2+\bar{\bar{\mathcal{A}}}_3\big)^\top\check{\Phi}+\bar{\mathcal{C}}_2^\top\check{\Lambda}
                  +\big(2\bar{\bar{\mathcal{F}}}_3^\top-\bar{\mathcal{F}}_3^\top\big)\check{v}_3+\bar{\bar{f}}_2\Big\}dt,\quad t\in[0,T],\\
   \check{X}_2(0)&=X_0,\quad \check{\Phi}_2(T)=0,
    \end{aligned}
\right.
\end{equation}
respectively. The problem of Player 2 could be solved in the following theorem.

\begin{theorem}
Let $P_1(\cdot), P_2(\cdot)$ satisfy (\ref{asymmRiccati21}) and (\ref{asymmRiccati22}), Player 2's optimal strategy $v_2^*(\cdot)$ is given by (\ref{asymmvbiaoda21}), where $(\hat{X}_2(\cdot), \check{X}_{2}(\cdot), \hat{\Phi}(\cdot), \hat{\Lambda}(\cdot))$ are the unique $\mathcal{G}_t^2$-adapted solutions to (\ref{asymmwenti21}), (\ref{asymmwenti22}).
\end{theorem}

\subsection{Problem of Player 3}

This section focuses on solving the problem of Player 3. Putting Player 2's optimal strategy $v_2^*(\cdot)$ into (\ref{asymmkuowei2}), the ``state" equation of Player 3 becomes:

\begin{equation}\label{asymmstate3}
\left\{
    \begin{aligned}
  dX_3(t)&=\Big\{\bar{\bar{\mathcal{A}}}_1X_3+\bar{\bar{\mathcal{A}}}_2\hat{X}_3+\bar{\bar{\mathcal{A}}}_3\check{X}_3+\bar{\mathcal{F}}_1\Phi
          +(\bar{\bar{\mathcal{F}}}_1-\bar{\mathcal{F}}_1)\hat{\Phi}+\bar{\mathcal{B}}_3v_3+\bar{\bar{b}}_2\Big\}dt\\
         &\quad +\sum_{i=1}^3 \big[\bar{\mathcal{C}}_iX_3+\bar{\sigma}_i\big]dW_i(t),\\
-d\Phi(t)&=\Big\{\mathcal{H}X_3-\mathcal{H}\hat{X}_3+\bar{\bar{\mathcal{A}}}_1^\top\Phi+\big(\bar{\bar{\mathcal{A}}}_2^\top-\bar{\mathcal{F}}_2^\top R_2^{-1}\bar{\mathcal{B}}_2^\top\big)\hat{\Phi}
          +\big(\bar{\bar{\mathcal{A}}}_3^\top+\bar{\mathcal{F}}_2^\top R_2^{-1}\bar{\mathcal{B}}_2^\top\big)\check{\Phi}\\
         &\qquad +\bar{\mathcal{C}}_2^\top\Lambda+\bar{\bar{\mathcal{F}}}_3^\top v_3
          +\big(\bar{\bar{\mathcal{F}}}_3-\bar{\mathcal{F}}_3\big)^\top\check{v}_3+\bar{\bar{f}}_2\Big\}dt-\Lambda(t)dW_2(t),\quad t\in[0,T],\\
   X_3(0)&=X_{0},\quad \Phi(T)=0,
    \end{aligned}
\right.
\end{equation}
where $X_3(\cdot)\triangleq X_2^{v_1^*,v_2^*,v_3}(\cdot)\equiv(x_2^*(\cdot),\psi_2(\cdot))^\top$, $\hat{X}_3(\cdot)\triangleq\hat{X}_2^{v_1^*,v_2^*,v_3}(\cdot)$, $\check{X}_3(\cdot)\triangleq\check{X}_2^{v_1^*,v_2^*,v_3}(\cdot)$.

\begin{remark}
Now, the dimension of $X_3(\cdot)$ is $2n\times 1$. Applying Lemma 5.4 in \cite{Xiong2008} to the state equation (\ref{asymmstate1}) and (\ref{asymmvarphi}), we can find that the equations about $(\check{X}_3(\cdot),\check{\Phi}(\cdot))$ and the equations about $(\hat{X}_3(\cdot),\hat{\Phi}(\cdot))$ are decoupled. According to (\ref{asymmvbiaoda21}), for any given $v_2(\cdot)$, the solution of (\ref{asymmstate3}) can be guaranteed though it is fully coupled.
\end{remark}

The cost functional of Player 3 can also be expanded accordingly. Setting
\begin{equation*}
\begin{split}
\bar{\mathcal{Q}}_3\triangleq\left(
     \begin{matrix}
     Q_3& 0\\
     0 & 0
     \end{matrix}
     \right),\quad
\bar{\mathcal{G}}_3\triangleq\left(
     \begin{matrix}
     G_{3}& 0\\
     0 & 0
     \end{matrix}
     \right),\quad
\bar{m}_3\triangleq\left(
     \begin{matrix}
     m_{3}\\ 0
     \end{matrix}
     \right),
\end{split}
\end{equation*}
we have
\begin{equation}
\begin{split}
 &J_3(v_1(\cdot)^*,v_2(\cdot)^*,v_3(\cdot))=\frac{1}{2}\mathbb{E}\bigg\{\int_0^T\Big[\big\langle\bar{\mathcal{Q}}_3(t)X_3(t),X_3(t)\big\rangle+\big\langle R_3(t)v_3(t),v_3(t)\big\rangle\\
 &\qquad +\big\langle 2\bar{m}_3(t),X_3(t)\big\rangle+\big\langle 2n_3(t),v_3(t)\big\rangle\Big]dt+\big\langle\bar{\mathcal{G}}_3X_3(T),X_3(T)\big\rangle\bigg\}.
\end{split}
\end{equation}

\begin{remark}
Different from the existing literature, there are two different state filtering in (\ref{asymmstate3}). We need to use the following equalities in the following derivation process of this subsection. For any $\mathcal{F}_t$-adapted processes $\xi(\cdot),\eta(\cdot)$, we have
\begin{equation}
\begin{aligned}
&\mathbb{E}\int_0^T\langle\mathbb{E}[\xi(t)|\mathcal{G}_t^1],\eta(t)\rangle dt=\mathbb{E}\int_0^T\langle\xi(t),\mathbb{E}[\eta(t)|\mathcal{G}_t^1]\rangle dt,\\
&\mathbb{E}\int_0^T\langle\mathbb{E}[\xi(t)|\mathcal{G}_t^2],\eta(t)\rangle dt=\mathbb{E}\int_0^T\langle\xi(t),\mathbb{E}[\eta(t)|\mathcal{G}_t^2]\rangle dt,\\
&\mathbb{E}\big[\mathbb{E}[\xi(t)|\mathcal{G}_t^2]|\mathcal{G}_t^1\big]=\mathbb{E}[\xi(t)|\mathcal{G}_t^1],\quad \mathbb{E}\big[\mathbb{E}[\xi(t)|\mathcal{G}_t^1]|\mathcal{G}_t^2]=\mathbb{E}[\xi(t)|\mathcal{G}_t^1\big].
\end{aligned}
\end{equation}
\end{remark}

Define the Hamiltionian function $H_3:[0,T]\times\mathbb{R}^{2n}\times\mathbb{R}^{2n}\times\mathbb{R}^{2n}\times\mathbb{R}^{2n}\times\mathbb{R}^{2n}\times\mathbb{R}^n\times\mathbb{R}^{2n}\times\mathbb{R}^{2n}\times\mathbb{R}^{2n}\times\mathbb{R}^{2n}
\times\mathbb{R}^{2n}\rightarrow\mathbb{R}$ of Player 3 as:
\begin{equation}
\begin{aligned}
  &H(t,X_3,\Phi,\hat{\Phi},\check{\Phi},\Lambda,v_3,Y_3,Z_1,Z_2,Z_3,\Psi_3)\\
  &\triangleq \Big[\bar{\bar{\mathcal{A}}}_1(t)X_3+\bar{\bar{\mathcal{A}}}_2(t)\hat{X}_3+\bar{\bar{\mathcal{A}}}_3(t)\check{X}_3+\bar{\mathcal{F}}_{1}(t)\Phi+(\bar{\bar{\mathcal{F}}}_1-\bar{\mathcal{F}}_1(t))\hat{\Phi}
   +\bar{\mathcal{B}}_3(t)v_3+\bar{\bar{b}}_2(t)\Big]^\top Y_3\\
  &\quad +\sum_{i=1}^3 \big[\bar{\mathcal{C}}_1(t)X_3+\bar{\sigma}_i(t)]^\top Z_i+\Big[\mathcal{H}(t)X_3-\mathcal{H}(t)\hat{X}_3+\bar{\bar{\mathcal{A}}}_1(t)^\top\Phi\\
  &\qquad +\big(\bar{\bar{\mathcal{A}}}_2(t)-\bar{\mathcal{F}}_2(t)^\top R_2^{-1}(t)\bar{\mathcal{B}}_2(t)^\top\big)\hat{\Phi}
   +\big(\bar{\bar{\mathcal{A}}}_3(t)+\bar{\mathcal{F}}_2(t)^\top R_2^{-1}(t)\bar{\mathcal{B}}_2(t)^\top\big)\check{\Phi}\\
  &\qquad +\bar{\mathcal{C}}_2(t)^\top\Lambda+\bar{\bar{\mathcal{F}}}_3(t)^\top v_3+\big(\bar{\bar{\mathcal{F}}}_3(t)^\top-\bar{\mathcal{F}}_3(t)^\top\big)\check{v}_3+\bar{\bar{f}}_2(t)\Big]^\top\Psi_3\\
  &\quad +\bigg[\frac{1}{2}\big\langle\bar{\mathcal{Q}}_3(t)X_3,X_3)\big\rangle+\frac{1}{2}\big\langle R_3(t)v_3,v_3\big\rangle+\big\langle\bar{m}_3(t),X_3\big\rangle+\big\langle n_3(t),v_3\big\rangle\bigg].
\end{aligned}
\end{equation}

The optimal strategy $v_3^*(\cdot)$ of Player 3 can be obtained by the maximum principle of FBSDEs with asymmetric information (see also, for example, \cite{SWX2016}, \cite{SWX2017}):
\begin{equation}
R_3v_3^*(t)+\bar{\mathcal{B}}_3^\top Y_3+\bar{\bar{\mathcal{F}}}_3\Psi_3+\big(\bar{\bar{\mathcal{F}}}_3-\bar{\mathcal{F}}_3\big)\check{\Psi}_3+n_3=0,\quad a.e.\ t\in[0,T],\ \ a.s.,
\end{equation}
where $(Y_3(\cdot), Z_1(\cdot), Z_2(\cdot), Z_3(\cdot), \Psi_3(\cdot))$ is the unique $\mathcal{F}_t$-adapted solution to the FBSDE:
\begin{equation}\label{asymmad3}
\left\{
    \begin{aligned}
        dY_3(t)=&-\Big[\bar{\bar{\mathcal{A}}}_1^\top Y_3+\bar{\bar{\mathcal{A}}}_2^\top\hat{Y}_3+\bar{\bar{\mathcal{A}}}_3^\top\check{Y}_3+\bar{\mathcal{F}}_1\Psi_3
                 +\big(\bar{\bar{\mathcal{F}}}_1-\bar{\mathcal{F}}_1\big)\hat{\Psi}_3+\sum_{i=1}^3 \bar{\mathcal{C}}_i(t)^\top Z_i\\
                &\quad +\bar{\mathcal{Q}}_3X_3^*+\tilde{m}_3\Big]dt+\sum_{i=1}^3 Z_i(t)dW_i(t),\\
     -d\Psi_3(t)&=-\Big[\bar{\bar{\mathcal{A}}}_1\Psi_3+\big(\bar{\bar{\mathcal{A}}}_2^\top-\bar{\mathcal{B}}_2R_2^{-1}\bar{\mathcal{F}}_2\big)\hat{\Psi}_3
                 +\big(\bar{\bar{\mathcal{A}}}_3^\top+\bar{\mathcal{B}}_2R_2^{-1}\bar{\mathcal{F}}_2\big)\check{\Psi}_3+\bar{\bar{\mathcal{F}}}_1Y_3\\
                &\qquad +\big(\bar{\bar{\mathcal{F}}}_1-\bar{\mathcal{F}}_1\big)\hat{Y}_3\Big]dt-\bar{\mathcal{C}}_2\Psi_3 dW_2(t),\quad t\in[0,T],\\
          Y_3(T)&=\bar{\mathcal{G}}_3X_3^*(T),\quad \Psi(0)=0,
    \end{aligned}
\right.
\end{equation}
where $(X_3^*(\cdot),\Phi^*(\cdot),\Lambda^*(\cdot))$ is the optimal state of Player 3.

Let $(X_3^*(\cdot),\Psi_3(\cdot))^\top$ be the optimal ``state" of Player 3, and set
\begin{equation}
    \begin{aligned}
    &\mathcal{X}_3\triangleq\left(\begin{matrix}X_3^*\\\Psi_3\end{matrix}\right),\
     \mathcal{Y}_3\triangleq\left(\begin{matrix}Y_{3}\\\Phi\end{matrix}\right),\
     \mathcal{Z}_1\triangleq\left(\begin{matrix}Z_1\\0\end{matrix}\right),\
     \mathcal{Z}_2\triangleq\left(\begin{matrix}Z_2\\\Lambda^*\end{matrix}\right),
     \mathcal{Z}_3\triangleq\left(\begin{matrix}Z_3\\0\end{matrix}\right),\
     \mathcal{X}_0=\left(\begin{matrix}X_0\\0\end{matrix}\right),\\
    &\bar{\bar{\mathfrak{A}}}_1\triangleq\left(\begin{matrix}\bar{\bar{\mathcal{A}}}_1& 0\\&\bar{\bar{\mathcal{A}}}_1\end{matrix}\right),\
     \bar{\bar{\mathfrak{A}}}_2\triangleq\left(\begin{matrix}\bar{\bar{\mathcal{A}}}_2& 0\\0 &\bar{\bar{\mathcal{A}}}_2-\bar{\mathcal{B}}_2R_2^{-1}\bar{\mathcal{F}}_2^\top\end{matrix}\right),\
     \bar{\bar{\mathfrak{A}}}_3\triangleq\left(\begin{matrix}\bar{\bar{\mathcal{A}}}_3& 0\\0 &\bar{\bar{\mathcal{A}}}_3+\bar{\mathcal{B}}_2R_2^{-1}\bar{\mathcal{F}}_2^\top\end{matrix}\right),\\
    &\bar{\mathfrak{F}}_1\triangleq\left(\begin{matrix}0& \bar{\mathcal{F}}_1\\\bar{\mathcal{F}}_1& 0\end{matrix}\right),\
     \bar{\bar{\mathfrak{F}}}_1\triangleq\left(\begin{matrix}0& \bar{\bar{\mathcal{F}}}_1\\\bar{\bar{\mathcal{F}}}_1& 0\end{matrix}\right),\
     \bar{\mathfrak{B}}_3\triangleq\left(\begin{matrix}\bar{\mathcal{B}}_3\\0\end{matrix}\right),\
     \bar{\mathfrak{C}}_1\triangleq\left(\begin{matrix}\bar{\mathcal{C}}_1& 0\\0& 0\end{matrix}\right),\
     \bar{\mathfrak{C}}_2\triangleq\left(\begin{matrix}\bar{\mathcal{C}}_2& 0\\0& \bar{\mathcal{C}}_2\end{matrix}\right),\\
    &\bar{\mathfrak{C}}_3\triangleq\left(\begin{matrix}\bar{\mathcal{C}}_3& 0\\0& 0\end{matrix}\right),\
     \bar{\mathfrak{Q}}_3\triangleq\left(\begin{matrix}\bar{\mathcal{Q}}_3& 0\\\mathcal{H}& 0\end{matrix}\right),\
     \bar{\bar{\mathfrak{Q}}}_3\triangleq\left(\begin{matrix}0\\\mathcal{H}\end{matrix}\right),\
     \bar{\mathfrak{G}}_3\triangleq\left(\begin{matrix}\bar{\mathcal{G}}_3& 0\\0& 0\end{matrix}\right),\
     \bar{\mathfrak{F}}_3\triangleq\left(\begin{matrix}0 & \bar{\mathcal{F}}_3\end{matrix}\right),\\
    &\bar{\bar{\mathfrak{F}}}_3\triangleq\left(\begin{matrix}0 & \bar{\bar{\mathcal{F}}}_3\end{matrix}\right),\
     \bar{\bar{b}}_3\triangleq\left(\begin{matrix}\bar{\bar{b}}_2\\0\end{matrix}\right),\
     \bar{\Sigma}_1\triangleq\left(\begin{matrix}\bar{\sigma}_1\\\Lambda^*\end{matrix}\right),\
     \bar{\Sigma}_2\triangleq\left(\begin{matrix}\bar{\sigma}_2\\0\end{matrix}\right),\
     \bar{\Sigma}_3\triangleq\left(\begin{matrix}\bar{\sigma}_3\\0\end{matrix}\right),\
     \bar{\bar{f}}_3\triangleq\left(\begin{matrix}\bar{m}_3\\\bar{\bar{f}}_2\end{matrix}\right).
    \end{aligned}
\end{equation}
Then (\ref{asymmstate3}) and (\ref{asymmad3}) can be rewritten as
\begin{equation}\label{asymmkuowei3}
\left\{
    \begin{aligned}
  d\mathcal{X}_3(t)&=\Big[\bar{\bar{\mathfrak{A}}}_1\mathcal{X}_3+\bar{\bar{\mathfrak{A}}}_2\hat{\mathcal{X}}_3+\bar{\bar{\mathfrak{A}}}_3\check{\mathcal{X}}_3+\bar{\mathfrak{F}}_1\mathcal{Y}_3
                    +\big(\bar{\bar{\mathfrak{F}}}_1-\bar{\mathfrak{F}}_1\big)\hat{\mathcal{Y}}_3+\bar{\mathfrak{B}}_3v_3^*+\bar{\bar{b}}_2\Big]dt\\
                   &\quad +\sum_{i=1}^3 \big[\bar{\mathfrak{C}}_i\mathcal{X}_3+\bar{\Sigma}_i\big]dW_i(t),\\
 -d\mathcal{Y}_3(t)&=\Big[\bar{\mathfrak{Q}}_3\mathcal{X}_3+\bar{\bar{\mathfrak{Q}}}_3\hat{\mathcal{X}}_3+\bar{\bar{\mathfrak{A}}}_1^\top\mathcal{Y}_3
                    +\bar{\bar{\mathfrak{A}}}_2^\top\hat{\mathcal{Y}}_3+\bar{\bar{\mathfrak{A}}}_3^\top\check{\mathcal{Y}}_3+\sum_{i=1}^3 \bar{\mathfrak{C}}_i^\top\mathcal{Z}_i+\bar{\bar{\mathfrak{F}}}_3^\top v_3^*\\
                   &\qquad +\big(\bar{\bar{\mathfrak{F}}}_3-\bar{\mathfrak{F}}_3)^\top\check{v}_3^*+\bar{\bar{f}}_3\Big]dt-\sum_{i=1}^3 \mathcal{Z}_i(t)dW_i(t),\quad t\in[0,T],\\
   \mathcal{X}_3(0)&=\mathcal{X}_0,\quad \mathcal{Y}_3(T)=\mathfrak{G}_3\mathcal{X}_3(T).
    \end{aligned}
\right.
\end{equation}
and similarly, the optimal strategy $v_3^*(\cdot)$ of Player 3 can be written as:
\begin{equation}\label{asymmv3chu}
  v_3^*(t)=-R_3^{-1}\big[\bar{\mathfrak{B}}_3^\top\mathcal{Y}_3+\bar{\bar{\mathfrak{F}}}_3\mathcal{X}_3+\big(\bar{\bar{\mathfrak{F}}}_3-\bar{\mathfrak{F}}_3\big)\check{\mathcal{X}}_3+n_3\big],\ a.e.\ t\in[0,T],\ \ a.s..
\end{equation}

Next, in order to represent the optimal strategy $v_3^*(\cdot)$ of Player 3 in the state feedback form, we assume that
\begin{equation}\label{asymmy3}
\mathcal{Y}_3(t)=\mathcal{P}_1(t)\mathcal{X}_3(t)+\mathcal{P}_2(t)\hat{\mathcal{X}}_3(t)+\mathcal{P}_3(t)\check{\mathcal{X}}_3(t)+\Omega(t),\quad t\in[0,T],
\end{equation}
where $\mathcal{P}_1(\cdot), \mathcal{P}_2(\cdot), \mathcal{P}_3(\cdot)$ are deterministic, differentiable functions with $\mathcal{P}_1(T)=\mathcal{G}_3$, $\mathcal{P}_2(T)=0$, $\mathcal{P}_3(T)=0$, and $\mathcal{F}_{t}$-adapted process pair $(\Omega(\cdot),\Pi(\cdot))$ satisfies the BSDE:
\begin{equation}
\left\{
    \begin{aligned}
       d\Omega(t)&=\Delta(t)dt+\Pi(t)dW_1(t),\\
        \Omega(T)&=0,
    \end{aligned}
\right.
\end{equation}
where $\Delta(\cdot)$ is some $\mathcal{F}_t$-adapted process to be determined.

Let (\ref{asymmy3}) be filtered on $\mathcal{G}_t^2$ to obtain $\hat{\mathcal{Y}}_3(\cdot)$ and on $\mathcal{G}_t^1$ to obtain $\check{\mathcal{Y}}_3(\cdot)$. Substituting $\mathcal{Y}_3(\cdot)$ into (\ref{asymmv3chu}), we can get the optimal state feedback strategy $v_3^*(\cdot)$ of Player 3 is:
\begin{equation}\label{asymmvbiaoda3}
\begin{split}
v_3^*(t)&=-R_3^{-1}\Big\{\big[\bar{\mathfrak{B}}_3^\top\mathcal{P}_1+\bar{\bar{\mathfrak{F}}}_3\big]\mathcal{X}_3+\bar{\mathfrak{B}}_3^\top P_2\hat{\mathcal{X}}_3+\big[\bar{\mathfrak{B}}_3^\top P_3\\
        &\qquad\qquad +\big(\bar{\bar{\mathfrak{F}}}_3-\bar{\mathfrak{F}}_3\big)\big]\check{\mathcal{X}}_3+\bar{\mathfrak{B}}_3^\top\Omega+n_3\Big\},\quad a.e.\ t\in[0,T],\ \ a.s..
\end{split}
\end{equation}

Using the same method as in sections 3.1 and 3.2, we introduce a system of three Riccati equations as follows:
\begin{equation}\label{asymmRiccati31}
\left\{
    \begin{aligned}
       &\dot{\mathcal{P}}_1+\mathcal{P}_1\big[\bar{\bar{\mathfrak{A}}}_1-\bar{\mathfrak{B}}_3R_3^{-1}\bar{\bar{\mathfrak{F}}}_3\big]
        +\big[\bar{\bar{\mathfrak{A}}}_1^\top-\bar{\bar{\mathfrak{F}}}_3^\top R_3^{-1}\bar{\mathfrak{B}}_3^\top\big]\mathcal{P}_1+\sum_{i=1}^3 \bar{\mathfrak{C}}_i^\top\mathcal{P}_1\bar{\mathfrak{C}}_i\\
       &\quad +\mathcal{P}_1\big[\bar{\mathfrak{F}}_1-\bar{\mathfrak{B}}_3R_3^{-1}\bar{\mathfrak{B}}_3^\top\big]\mathcal{P}_1+\bar{\mathfrak{Q}}_3-\bar{\bar{\mathfrak{F}}}_3^\top R_3^{-1}\bar{\bar{\mathfrak{F}}}_3=0,\\
       &\mathcal{P}_1(T)=\bar{\mathfrak{G}}_3,
    \end{aligned}
\right.
\end{equation}
\begin{equation}\label{asymmRiccati32}
\left\{
    \begin{aligned}
       &\dot{\mathcal{P}}_2+\mathcal{P}_2\big[\bar{\bar{\mathfrak{A}}}_1+\bar{\bar{\mathfrak{A}}}_2-\bar{\mathfrak{B}}_3R_3^{-1}\bar{\bar{\mathfrak{F}}}_3\big]
        +\big[\bar{\bar{\mathfrak{A}}}_1^\top+\bar{\bar{\mathfrak{A}}}_2^\top-\bar{\bar{\mathfrak{F}}}_3^\top R_3^{-1}\bar{\mathfrak{B}}_3^\top\big]\mathcal{P}_2\\
       &\quad +\mathcal{P}_1\big[\bar{\bar{\mathfrak{F}}}_1-\bar{\mathfrak{B}}_3R_3^{-1}\bar{\mathfrak{B}}_3^\top\big]\mathcal{P}_2
        +\mathcal{P}_2\big[\bar{\bar{\mathfrak{F}}}_1-\bar{\mathfrak{B}}_3R_3^{-1}\bar{\mathfrak{B}}_3^\top\big]\mathcal{P}_1+\sum_{i=2}^3 \bar{\mathfrak{C}}_i^\top\mathcal{P}_2\bar{\mathfrak{C}}_i\\
       &\quad +\mathcal{P}_2\big[\bar{\bar{\mathfrak{F}}}_1-\bar{\mathfrak{B}}_3R_3^{-1}\bar{\mathfrak{B}}_3^\top\big]\mathcal{P}_2+\bar{\bar{\mathfrak{Q}}}_3+\mathcal{P}_1\bar{\bar{\mathfrak{A}}}_2
        +\bar{\bar{\mathfrak{A}}}_2^\top\mathcal{P}_1+\mathcal{P}_1\big(\bar{\bar{\mathfrak{F}}}_1-\bar{\mathfrak{F}}_1\big)\mathcal{P}_1=0,\\
       &\mathcal{P}_2(T)=0,
    \end{aligned}
\right.
\end{equation}
\begin{equation}\label{asymmRiccati33}
\left\{
    \begin{aligned}
       &\dot{\mathcal{P}}_3+\mathcal{P}_3\big[\bar{\bar{\mathfrak{A}}}_1+\bar{\bar{\mathfrak{A}}}_2+\bar{\bar{\mathfrak{A}}}_3
        -\bar{\mathfrak{B}}_3R_3^{-1}\big(2\bar{\bar{\mathfrak{F}}}_3-\bar{\mathfrak{F}}_3\big)\big]
        +\big[\bar{\bar{\mathfrak{A}}}_1^\top+\bar{\bar{\mathfrak{A}}}_2^\top+\bar{\bar{\mathfrak{A}}}_3^\top-\big(2\bar{\bar{\mathfrak{F}}}_3
        -\bar{\mathfrak{F}}_3\big)^\top R_3^{-1}\bar{\mathfrak{B}}_3^\top\big]\mathcal{P}_3\\
       &\quad +\mathcal{P}_1\big[\bar{\bar{\mathfrak{A}}}_3-\bar{\mathfrak{B}}_3R_3^{-1}\big(\bar{\bar{\mathfrak{F}}}_3-\bar{\mathfrak{F}}_3\big)\big]
        +\big[\bar{\bar{\mathfrak{A}}}_3^\top-\big(\bar{\bar{\mathfrak{F}}}_3-\bar{\mathfrak{F}}_3\big)^\top R_3^{-1}\bar{\mathfrak{B}}_3^\top\big]\mathcal{P}_1\\
       &\quad +\mathcal{P}_1\big[\bar{\bar{\mathfrak{F}}}_1-\bar{\mathfrak{B}}_3R_3^{-1}\bar{\mathfrak{B}}_3^\top\big]\mathcal{P}_3
        +\mathcal{P}_3\big[\bar{\bar{\mathfrak{F}}}_1-\bar{\mathfrak{B}}_3R_3^{-1}\bar{\mathfrak{B}}_3^\top\big]\mathcal{P}_3\\
       &\quad +\mathcal{P}_2\big[\bar{\bar{\mathfrak{A}}}_3-\bar{\mathfrak{B}}_3R_3^{-1}\big(\bar{\bar{\mathfrak{F}}}_3-\bar{\mathfrak{F}}_3\big)\big]
        +\big[\bar{\bar{\mathfrak{A}}}_3^\top-\big(\bar{\bar{\mathfrak{F}}}_3-\bar{\mathfrak{F}}_3\big)^\top R_3^{-1}\bar{\mathfrak{B}}_3^\top\big]\mathcal{P}_{2}\\
       &\quad +\mathcal{P}_2\big[\bar{\bar{\mathfrak{F}}}_1-\bar{\mathfrak{B}}_3R_3^{-1}\bar{\mathfrak{B}}_3^\top\big]\mathcal{P}_3
        +\mathcal{P}_3\big[\bar{\bar{\mathfrak{F}}}_1-\bar{\mathfrak{B}}_3R_3^{-1}\bar{\mathfrak{B}}_3^\top\big]\mathcal{P}_2+\bar{\mathfrak{C}}_3^\top\mathcal{P}_3\bar{\mathfrak{C}}_3\\
       &\quad +\mathcal{P}_3\big[\bar{\bar{\mathfrak{F}}}_1-\bar{\mathfrak{B}}_3R_3^{-1}\bar{\mathfrak{B}}_3^\top\big]\mathcal{P}_3
        +(\bar{\bar{\mathfrak{F}}}_3^\top-\bar{\mathfrak{F}}_3^\top\big]R_3^{-1}\bar{\bar{\mathfrak{F}}}_3-\bar{\bar{\mathfrak{F}}}_3^\top R_3^{-1}\big[\bar{\bar{\mathfrak{F}}}_3-\bar{\mathfrak{F}}_3\big]=0,\\
       &\mathcal{P}_3(T)=0.
    \end{aligned}
\right.
\end{equation}
It can be easily seen that the system (\ref{asymmRiccati31})-(\ref{asymmRiccati33}) admits a differentiable solution triple $(\mathcal{P}_1(\cdot)$, $\mathcal{P}_2(\cdot)$, $\mathcal{P}_3(\cdot))$, since we can solve them sequently. Then, we can get $(\Omega(\cdot),\Pi(\cdot))$ satisfies the following BSDE:
\begin{equation}\label{asymmOmega31}
\left\{
    \begin{aligned}
      -d\Omega(t)&=\bigg\{\big[\bar{\bar{\mathfrak{A}}}_1^\top-\bar{\bar{\mathfrak{F}}}_3^\top R_3^{-1}(t)\bar{\mathfrak{B}}_3^\top
                  +\mathcal{P}_1\big(\bar{\mathfrak{F}}_1-\bar{\mathfrak{B}}_3R_3^{-1}\bar{\mathfrak{B}}_3^\top\big)\big]\Omega\\
                 &\qquad +\big[\bar{\bar{\mathfrak{A}}}_2^\top+\mathcal{P}_1\big(\bar{\bar{\mathfrak{F}}}_1-\bar{\mathfrak{F}}_1\big)
                  +\mathcal{P}_2\big(\bar{\bar{\mathfrak{F}}}_1-\bar{\mathfrak{B}}_3R_3^{-1}\bar{\mathfrak{B}}_3^\top\big)\big]\hat{\Omega}\\
                 &\qquad +\big[\bar{\bar{\mathfrak{A}}}_3^\top-\big(\bar{\bar{\mathfrak{F}}}_3-\bar{\mathfrak{F}}_3\big)^\top R_3^{-1}\bar{\mathfrak{B}}_3^\top
                  +\mathcal{P}_3\big(\bar{\bar{\mathfrak{F}}}_1-\bar{\mathfrak{B}}_3R_3^{-1}\bar{\mathcal{B}}_3^\top\big)\big]\check{\Omega}+\bar{\mathfrak{C}}_2^\top\Pi\\
                 &\qquad +\sum_{i=1}^3 \bar{\mathfrak{C}}_i^\top\mathcal{P}_1\bar{\Sigma}_i+\sum_{i=2}^3 \bar{\mathfrak{C}}_i^\top\mathcal{P}_2\bar{\Sigma}_i
                  +\bar{\mathfrak{C}}_3^\top\mathcal{P}_3\bar{\Sigma}_3+\mathcal{P}_1\bar{\bar{b}}_2+\bar{\bar{f}}_3\\
                 &\qquad -\bigg[\sum_{i=1}^3 \mathcal{P}_i\bar{\mathfrak{B}}_3+\bar{\mathfrak{F}}_3^\top\bigg]R_3^{-1}n_3\bigg\}dt-\Pi(t)dW(t),\quad t\in[0,T],\\
        \Omega(T)&=0.
    \end{aligned}
\right.
\end{equation}

\begin{remark}
Similarly, readers can refer to \cite{SWX2017} for more details about the derivation of (\ref{asymmRiccati31})-(\ref{asymmRiccati33}). In this way, it can also be shown the solvability of (\ref{asymmkuowei3}), since  $(\mathcal{X}_3(\cdot)$, $\hat{\mathcal{X}}_3(\cdot)$, $\check{\mathcal{X}}_3(\cdot)$, $\Omega(\cdot)$, $\Pi(\cdot))$ in (\ref{asymmvbiaoda3}) are the unique solutions to the following equations (\ref{asymmwenti31})-(\ref{asymmwenti33}).
\end{remark}

Applying Lemma 5.4 in \cite{Xiong2008} to state equation (\ref{asymmkuowei3}) and (\ref{asymmOmega31}), we derive three optimal filtering equations:
\begin{equation}\label{asymmwenti31}
\left\{
    \begin{aligned}
   d\mathcal{X}_3(t)&=\Big\{\big[\bar{\bar{\mathfrak{A}}}_1-\bar{\mathfrak{B}}_3R_3^{-1}\bar{\bar{\mathfrak{F}}}_3
                     +\big(\bar{\mathfrak{F}}_1-\bar{\mathfrak{B}}_3R_3^{-1}\bar{\mathfrak{B}}_3^\top\big)\mathcal{P}_1\big]\mathcal{X}_3\\
                    &\qquad +\big[\bar{\bar{\mathfrak{A}}}_2+\big(\bar{\bar{\mathfrak{F}}}_1-\bar{\mathfrak{F}}_1\big)\mathcal{P}_1
                     +\big(\bar{\bar{\mathfrak{F}}}_1-\bar{\mathfrak{B}}_3R_3^{-1}\bar{\mathfrak{B}}_3^\top\big)\mathcal{P}_2\big]\hat{\mathcal{X}}_3\\
                    &\qquad +\big[\bar{\bar{\mathfrak{A}}}_3-\bar{\mathfrak{B}}_3R_3^{-1}\big(\bar{\bar{\mathfrak{F}}}_3-\bar{\mathfrak{F}}_3\big)
                     +\big(\bar{\bar{\mathfrak{F}}}_1-\bar{\mathfrak{B}}_3R_3^{-1}\bar{\mathfrak{B}}_3^\top\big)\mathcal{P}_3\big]\check{\mathcal{X}}_3\\
                    &\qquad +\big[\bar{\mathfrak{F}}_1-\bar{\mathfrak{B}}_3R_3^{-1}\bar{\mathfrak{B}}_3^\top\big]\Omega+\big[\bar{\bar{\mathfrak{F}}}_1-\bar{\mathfrak{F}}_1\big]\hat{\Omega}
                     +\bar{\bar{b}}_2-\bar{\mathfrak{B}}_3R_3^{-1}n_3\Big\}dt\\
                    &\quad +\sum_{i=1}^3 \big[\bar{\mathfrak{C}}_i\mathcal{X}_3+\bar{\Sigma}_i\big]dW_i(t),\\
         -d\Omega(t)&=\bigg\{\big[\bar{\bar{\mathfrak{A}}}_1^\top-\bar{\bar{\mathfrak{F}}}_3^\top R_3^{-1}\bar{\mathfrak{B}}_3^\top
                     +\mathcal{P}_1\big(\bar{\mathfrak{F}}_1-\bar{\mathfrak{B}}_3R_3^{-1}\bar{\mathfrak{B}}_3^\top\big)\big]\Omega\\
                    &\qquad +\big[\bar{\bar{\mathfrak{A}}}_2^\top+\mathcal{P}_1\big(\bar{\bar{\mathfrak{F}}}_1-\bar{\mathfrak{F}}_1\big)
                     +\mathcal{P}_2\big(\bar{\bar{\mathfrak{F}}}_1-\bar{\mathfrak{B}}_3R_3^{-1}\bar{\mathfrak{B}}_3^\top\big)\big]\hat{\Omega}\\
                    &\qquad +\big[\bar{\bar{\mathfrak{A}}}_3^\top-\big(\bar{\bar{\mathfrak{F}}}_3-\bar{\mathfrak{F}}_3\big)^\top R_3^{-1}\bar{\mathfrak{B}}_3^\top
                     +\mathcal{P}_3\big(\bar{\bar{\mathfrak{F}}}_1-\bar{\mathfrak{B}}_3R_3^{-1}\bar{\mathcal{B}}_3^\top\big)\big]\check{\Omega}\\
                    &\qquad +\sum_{i=1}^3 \bar{\mathfrak{C}}_i^\top\mathcal{P}_1\bar{\Sigma}_i+\sum_{i=2}^3 \bar{\mathfrak{C}}_i^\top\mathcal{P}_2\bar{\Sigma}_i
                     +\bar{\mathfrak{C}}_3^\top\mathcal{P}_3\bar{\Sigma}_3+\mathcal{P}_1\bar{\bar{b}}_2+\bar{\bar{f}}_3\\
                    &\qquad -\bigg[\sum_{i=1}^3 \mathcal{P}_i\bar{\mathfrak{B}}_3+\bar{\mathfrak{F}}_3^\top\bigg]R_3^{-1}n_3\bigg\}dt-\Pi(t)dW(t),\quad t\in[0,T],\\
    \mathcal{X}_3(0)&=\mathcal{X}_0,\quad \Omega(T)=0,
    \end{aligned}
\right.
\end{equation}
\begin{equation}\label{asymmwenti32}
\left\{
    \begin{aligned}
d\hat{\mathcal{X}}_3(t)&=\Big\{\big[\bar{\bar{\mathfrak{A}}}_1+\bar{\bar{\mathfrak{A}}}_2-\bar{\mathfrak{B}}_3R_3^{-1}\bar{\bar{\mathfrak{F}}}_3
                        +\big(\bar{\bar{\mathfrak{F}}}_1-\bar{\mathfrak{B}}_3R_3^{-1}\bar{\mathfrak{B}}_3^\top\big)\big(\mathcal{P}_1+\mathcal{P}_2\big)\big]\hat{\mathcal{X}}_3\\
                       &\qquad +\big[\bar{\bar{\mathfrak{A}}}_3-\bar{\mathfrak{B}}_3R_3^{-1}\big(\bar{\bar{\mathfrak{F}}}_3-\bar{\mathfrak{F}}_3\big)
                        +\big(\bar{\bar{\mathfrak{F}}}_1-\bar{\mathfrak{B}}_3R_3^{-1}\bar{\mathfrak{B}}_3^\top\big)\mathcal{P}_3\big]\check{\mathcal{X}}_3\\
                       &\qquad +\big[\bar{\bar{\mathfrak{F}}}_1-\bar{\mathfrak{B}}_3R_3^{-1}\bar{\mathfrak{B}}_3^\top\big]\hat{\Omega}-\bar{\mathfrak{B}}_3R_3^{-1}n_3+\bar{\bar{b}}_2\Big\}dt
                        +\sum_{i=2}^3 \big[\bar{\mathfrak{C}}_i\hat{\mathcal{X}}_3+\bar{\Sigma}_i\big]dW_i(t),\\
      -d\hat{\Omega}(t)&=\bigg\{\big[\bar{\bar{\mathfrak{A}}}_1^\top+\bar{\bar{\mathfrak{A}}}_2^\top-\bar{\bar{\mathfrak{F}}}_3^\top R_3^{-1}\bar{\mathfrak{B}}_3^\top
                        +\big(\mathcal{P}_1+\mathcal{P}_2\big)\big(\bar{\bar{\mathfrak{F}}}_1-\bar{\mathfrak{B}}_3R_3^{-1}\bar{\mathfrak{B}}_3^\top\big)\big]\hat{\Omega}\\
                       &\qquad +\big[\bar{\bar{\mathfrak{A}}}_3^\top-\big(\bar{\bar{\mathfrak{F}}}_3-\bar{\mathfrak{F}}_3\big)^\top R_3^{-1}\bar{\mathfrak{B}}_3^\top
                        +\mathcal{P}_3\big(\bar{\bar{\mathfrak{F}}}_1-\bar{\mathfrak{B}}_3R_3^{-1}\bar{\mathcal{B}}_3^\top\big)\big]\check{\Omega}+\bar{\mathfrak{C}}_2^\top\hat{\Pi}\\
                       &\qquad +\sum_{i=1}^3 \bar{\mathfrak{C}}_i^\top\big(\mathcal{P}_1+\mathcal{P}_2\big)\bar{\Sigma}_i+\mathcal{P}_1\bar{\bar{b}}_2
                        +\bar{\bar{f}}_3-\bigg[\sum_{i=1}^3 \mathcal{P}_i\bar{\mathfrak{B}}_3+\bar{\mathfrak{F}}_3^\top\bigg]R_3^{-1}n_3\bigg\}dt,\quad t\in[0,T],\\
 \hat{\mathcal{X}}_3(0)&=\mathcal{X}_0,\quad \hat{\Omega}(T)=0,
    \end{aligned}
\right.
\end{equation}
and
\begin{equation}\label{asymmwenti33}
\left\{
    \begin{aligned}
d\check{\mathcal{X}}_3(t)&=\bigg\{\bigg[\sum_{i=1}^3 \bar{\bar{\mathfrak{A}}}_i
                          -\bar{\mathfrak{B}}_3R_3^{-1}\big(2\bar{\bar{\mathfrak{F}}}_3-\bar{\mathfrak{F}}_3\big)
                          +\big(\bar{\bar{\mathfrak{F}}}_1-\bar{\mathfrak{B}}_3R_3^{-1}\bar{\mathfrak{B}}_3^\top\big)\bigg(\sum_{i=1}^3 \mathcal{P}_i\bigg)\bigg]\check{\mathcal{X}}_3\\
                         &\qquad +\big[\bar{\bar{\mathfrak{F}}}_1-\bar{\mathfrak{B}}_3R_3^{-1}\bar{\mathfrak{B}}_3^\top\big]\check{\Omega}-\bar{\mathfrak{B}}_3R_3^{-1}n_3+\bar{\bar{b}}_2\bigg\}dt
                          +\big[\bar{\mathfrak{C}}_3\check{\mathcal{X}}_3+\bar{\Sigma}_3]dW_3(t),\\
  -d\check{\Omega}(t)&=\bigg\{\bigg[\sum_{i=1}^3 \bar{\bar{\mathfrak{A}}}_i^\top
                          -\big(2\bar{\bar{\mathfrak{F}}}_3-\bar{\mathfrak{F}}_3\big)^\top R_3^{-1}\bar{\mathfrak{B}}_3^\top
                          +\bigg(\sum_{i=1}^3 \mathcal{P}_i\bigg)\big(\bar{\bar{\mathfrak{F}}}_1-\bar{\mathfrak{B}}_3R_3^{-1}\bar{\mathfrak{B}}_3^\top\big)\bigg]\check{\Omega}\\
                         &\qquad +\bar{\mathfrak{C}}_2^\top\check{\Pi}+\sum_{i=1}^3 \bar{\mathfrak{C}}_i^\top\big(\mathcal{P}_1+\mathcal{P}_2\big)\bar{\Sigma}_i+\mathcal{P}_1\bar{\bar{b}}_2+\bar{\bar{f}}_3\\
                         &\qquad -\bigg[\sum_{i=1}^3 \mathcal{P}_i\bar{\mathfrak{B}}_3+\bar{\mathfrak{F}}_3^\top\bigg]R_3^{-1}n_3\bigg\}dt,\quad t\in[0,T],\\
        \check{\mathcal{X}}_3(0)&=\mathcal{X}_0,\quad \check{\Omega}(T)=0.
    \end{aligned}
\right.
\end{equation}

The problem of Player 3 can be solved in the following theorem.
\begin{theorem}
Let $\mathcal{P}_1(\cdot)$, $\mathcal{P}_2(\cdot)$, $\mathcal{P}_3(\cdot)$ satisfy the system of Riccati equations (\ref{asymmRiccati31})-(\ref{asymmRiccati33}), Player 3's optimal strategy is given by (\ref{asymmvbiaoda3}), where $(\mathcal{X}_3(\cdot)$, $\hat{\mathcal{X}}_3(\cdot)$, $\check{\mathcal{X}}_3(\cdot)$, $\Omega(\cdot)$, $\Pi(\cdot))$ are the unique adapted solutions to (\ref{asymmwenti31})-(\ref{asymmwenti33}).
\end{theorem}

In addition, we can rewrite Player 1 and Player 2's optimal feedback strategies (\ref{asymmvbiaoda1}) and (\ref{asymmvbiaoda21}) as follows:
\begin{equation}\label{asymmv2zhong}
\begin{aligned}
v_2^*(t)&=-R_2^{-1}\bigg\{\left\{\begin{matrix}\left[\begin{matrix}\left(\begin{matrix}B_2&0\end{matrix}\right)P_1&0\end{matrix}\right]
                             +\left[\begin{matrix}0&\left(\begin{matrix}B_2&0\end{matrix}\right)\end{matrix}\right](\mathcal{P}_1+\mathcal{P}_2)\end{matrix}\right\}\hat{\mathcal{X}}_3(t)\\
        &\qquad\qquad\qquad +\left\{\begin{matrix}\left[\begin{matrix}\left(\begin{matrix}B_2&0\end{matrix}\right)P_2+\bar{\mathcal{F}}_2&0\end{matrix}\right]
         +\left[\begin{matrix}0&\left(\begin{matrix}B_2&0\end{matrix}\right)\end{matrix}\right]\mathcal{P}_3\end{matrix}\right\}\check{\mathcal{X}}_3(t)\\
        &\qquad\qquad\qquad +\left[\begin{matrix}0&\left(\begin{matrix}B_2&0\end{matrix}\right)\end{matrix}\right]\hat{\Omega}(t)+n_2(t)\bigg\},\quad a.e.\ t\in[0,T],\ \ a.s.,
\end{aligned}
\end{equation}
\begin{equation}\label{asymmv1zhong}
\begin{aligned}
v_1^*(t)&=-R_1^{-1}(t)\bigg\{\Big\{\left[\begin{matrix}\left(\begin{matrix}B_1^\top p_1&0\end{matrix}\right)
         +\left(\begin{matrix}0&B_1^\top\end{matrix}\right)P_2&0\end{matrix}\right]\\
        &\qquad\qquad\qquad +\left[\begin{matrix}0&\left(\begin{matrix}0&B_1^\top\end{matrix}\right)\end{matrix}\right](\mathcal{P}_1+\mathcal{P}_2+\mathcal{P}_3)\Big\}\check{\mathcal{X}}_3(t)\\
        &\qquad\qquad\qquad +\left[\begin{matrix}0&\left(\begin{matrix}0&B_1^\top\end{matrix}\right)\end{matrix}\right]\check{\Omega}(t)+n_1(t)\bigg\},\quad a.e.\ t\in[0,T],\ \ a.s..
\end{aligned}
\end{equation}
Up to now, we have obtained the optimal feedback equilibrium strategies of our stochastic LQ Stakelberg differential game with asymmetric information.

\section{Conslusions and Prospects}

In this paper, we have studied a three-level stochastic LQ Stackelberg differential game with asymmetric information, and the state feedback representation of the players' equilibrium strategies is given. The solution is performed to solve the optimal control problems of Player 1, Player 2 and Player 3 in turn. Maximum principle with partial information and optimal filtering are used to seek the equilibrium strategies of Player 1 and Player 2. When Player 1 and Player 2 exercise their optimal strategies, the state equation for the optimal control problem faced by Player 3 is a fully coupled FBSDE with two different types of stochastic filtering terms, which is different from existing literatures. A new system of high-dimensional Riccati equations is introduced to get Player 3's optimal state feedback strategy.

In the future, it is worthy to study the three-echelon supply chain management problem with asymmetric information, such as cooperative advertising and pricing problem in a dynamic three-echelon supply chain with asymmetric information.

\end{document}